\documentclass[11pt]{article}%revised version 2019-2-22
\def\q{\hfill\rule{1ex}{1ex}}
\def\0{\emptyset}

\def\q{\hfill\rule{1ex}{1ex}}

\def\n{\noindent}

\usepackage{graphicx}
\usepackage{authblk}
\usepackage{amsfonts}
\usepackage{enumitem}
\usepackage{calc}
\usepackage{amsthm}
\usepackage{amssymb}
\usepackage{amsmath}
\usepackage{amsthm}
\usepackage{tikz}
\usepackage{pgfplots}

\begin{document}
\title{\bf On the anti-Ramsey number of forests}
%\thanks{This work is partially supported by NNSFC (No. 11771247).}}
\author[1,2,4]{{\small\bf Chunqiu Fang}\thanks{ fcq15@mails.tsinghua.edu.cn, supported in part by CSC(No. 201806210164), and NNSFC (No. 11771247).}}
\author[2]{{\small\bf Ervin Gy\H ori}\thanks{ gyori.ervin@renyi.mta.hu, supported in part by the National Research, Development and Innovation Office NKFIH, grants K116769, K117879 and K126853.}}
\author[1]{{\small\bf Mei Lu}\thanks{ mlu@math.tsinghua.edu.cn,  supported in part by NNSFC (No. 11771247).}}
\author[2,3,5]{{\small\bf Jimeng Xiao}\thanks{ xiaojimeng@mail.nwpu.edu.cn, supported in part by CSC(No. 201706290171).}}

\affil[1]{Department of Mathematical Sciences, Tsinghua University, Beijing 100084, China }
\affil[2]{Alfr\'ed R\'enyi Institute of Mathematics, Hungarian Academy of Sciences, Re\'altanoda u.13-15, 1053 Budapest, Hungary}
\affil[3]{Department of Applied Mathematics, Northwestern Polytechnical University, Xi'an, China}
\affil[4]{Yau Mathematical Sciences Center, Tsinghua University, Beijing 10084, China}
\affil[5]{Xi'an-Budapest Joint Research Center for Combinatorics, Northwestern Polytechnical University, Xi'an, China}

\date{}

\maketitle\baselineskip 16.7pt

\begin{abstract}
We call a subgraph of an edge-colored graph rainbow subgraph, if all of its edges have different colors. The anti-Ramsey number of a graph $G$ in a complete graph $K_{n}$, denoted by $ar(K_{n}, G)$, is the maximum number of colors in an edge-coloring of $K_{n}$ with no rainbow subgraph copy of $G$. In this paper, we  determine the exact value of the anti-Ramsey number for star forests and the approximate value of the anti-Ramsey number for linear forests. Furthermore, we compute the exact value of
$ar(K_{n}, 2P_{4})$ for $n\ge 8$ and $ar(K_{n}, S_{p,q})$ for large $n$, where $S_{p,q}$ is the double star with $p+q$ leaves.

\end{abstract}

\vskip.2cm

\n{\large\bf 1. Introduction}

\vskip.2cm

Let $G $ be a simple undirected graph. For $x\in V (G)$, we denote the {\em neighborhood} (the set of neighbors of $x$) and the {\em degree} of $x$ in $G$ by $N_G(x)$ and $d_G(x)$, respectively.  The {\em maximum degree}  and the {\em minimum degree} of $G$ are denoted by $\Delta(G )$ and $\delta(G )$, respectively. %The {\em common neighborhood} of $U\subset V(G)$ is the set of vertices in $V(G)\setminus U$ that are adjacent to every vertex in $U$.
For $\emptyset\not= X \subset V (G)$, $G[X]$ is the subgraph of $G$ induced by $X$ and $G -X$ is the subgraph of $G$ induced by $V (G)\setminus X$. If $X=\{x\}$, then $G-X$ will be denoted by $G-x$ for short. Given a graph $G= (V, E)$, for any (not necessarily disjoint) vertex sets $A, B\subset V$,  let $E_{G}(A, B):=\{uv\in E(G)| u\ne v, u\in A, v\in B\}$.   A {\em star forest} is a forest whose  components are stars and a {\em linear forest} is a forest whose  components are paths. We use $tG$ and $\overline{G}$ to denote $t$ vertex-disjoint copies and the complement of $G$, respectively. Given two vertex disjoint graphs $G_{1}$ and $G_{2}$, we denote by $G_{1}+G_{2}$  the {\em join} of graphs $G_{1}$ and $G_{2}$, that is the graph obtained from $G_{1}\cup G_{2}$ by joining each vertex of $G_{1}$ with each vertex of $G_{2}$.

 We call a subgraph of an edge-colored graph {\em rainbow}, if all of its edges have different colors.  Let $G$ be a  graph. The {\em anti-Ramsey number} $ar(K_{n}, G)$ is the maximum number of colors in an edge-coloring of $K_{n}$ which has no rainbow copy of  $G$.  The {\em Tur\'an number} $ex(n, G)$ is the maximum number of edges of a simple graph on $n$ vertices without a copy of $G$.
The anti-Ramsey number was first studied by Erd\H os, Simonovits and S\'os \cite{Erdos}. They showed that the anti-Ramsey number is  closely related to Tur\'an number. Since then, there are plentiful results in this field, including  cycles \cite{Alon, Montellano-2}, cliques \cite{Montellano-3,Schiermeyer}, trees \cite{Jiang, Jiang-1}  and so on. See  Fujita, Magnant and Ozeki \cite{Fujita, FMO} for an abundant survey. Among these results, almost the considered graphs are connected graphs, and  a few unconnected graphs are considered including matchings \cite{Chen, Haas}, vertex-disjoint cliques \cite{Yuan}. In this paper, we will consider the cases that $G$ is a star forest or a linear forest.

We mention some of the results, which are relevant to our work.

Jiang \cite{Jiang} and Montellano-Ballesteros \cite{Montellano} independently found the anti-Ramsey number for stars.

\noindent{\bf Theorem 1.} (\cite{Jiang},\cite{Montellano}) {\em For $n\ge p+2\ge 3$, $$\left\lfloor\frac{(p-1)n}{2}\right\rfloor+\left\lfloor\frac{n}{n-p+1}\right\rfloor\le ar(K_{n}, K_{1,p+1})\le\left\lfloor\frac{(p-1)n}{2}+\frac{1}{2}\left\lfloor\frac{2n}{n-p+1}\right\rfloor\right\rfloor.$$}

By Theorem 1, $ar(K_{n}, K_{1,p})=\lfloor\frac{(p-2)n}{2}\rfloor+1$ for $n\ge 3p+4$ and $p\ge 2$.

Simonovits and S\'os \cite{Simonovits} considered the anti-Ramsey number for paths and obtained the following result.

\noindent{\bf Theorem 2.} (\cite{Simonovits}) {\em Let $P_{k+1}$ be a path of length $k\ge 2$ and $k\equiv r$ $ (\mod2), 0\le r\le 1$. For large enough $n$ $(n\ge \frac{5}{4}k+C$ for some universal constant $C)$, $$ar(K_{n}, P_{k+1})=\left(\left\lfloor \frac{k}{2}\right \rfloor-1\right)n-\binom{\lfloor \frac{k}{2} \rfloor}{2}+1+r.$$}

Let $\Omega_{k}$ denote the family of graphs that contain $k$ vertex-disjoint cycles. Jin and Li \cite{Jin} computed the anti-Ramsey number for $\Omega_{2}.$

\noindent{\bf Theorem 3.} (\cite{Jin}) {\em For any $n\ge 6$, $ar(K_{n}, \Omega_{2})=\max\{2n-2, 11\}.$}

The Tur\'an number of $tK_{2}$ was determined by Erd\H os and Gallai \cite{Erdos1} as $ex(n,tK_{2})=\max\{\binom{2t-1}{2}, \\ \binom{t-1}{2}+(t-1)(n-t+1)\}$ for $n\ge 2t\ge 2$.
The anti-Ramsey number of matchings was first considered by Schiermeyer \cite{Schiermeyer}.

\noindent{\bf Theorem 4.}(\cite{Schiermeyer}) {\em For $n\ge3t+3\ge 8$, we have $$ar(K_{n},tK_{2})=ex(n, (t-1)K_{2})+1=(t-2)n-\binom{t-1}{2}+1.$$}

 Later, Chen, Li and Tu \cite{Chen} and independently Fujita, Kaneko, Schiermeyer and Suzuki \cite{Fujita} showed that $ar(K_{n},tK_{2})=ex(n, (t-1)K_{2})+1$ for  $n\ge 2t+1\ge 5$. The values $ar(K_{2t}, tK_{2})=ex(2t, (t-1)K_{2})+1$ for $2\le t\le 6$ and $ar(K_{2t}, tK_{2})=ex(2t, (t-1)K_{2})+2$ for $t\ge 7$ were determined in \cite{Chen} and by Haas and Young \cite{Haas}, independently.

Gilboa and Roditty \cite{Gilboa} considered the graphs with small connected components and proved inductive results of the form ``if $ar(K_{n}, G\cup t_{0}P_{s})\le f(n, t_{0}, G)$ for sufficiently large $n$, then $ar(K_{n}, G\cup tP_{s})\le f(n, t, G)$ for sufficiently large $n$ and $t\ge t_{0}$, where $s=2$ or $3$."  These results imply the following theorem.

\noindent{\bf Theorem 5.} (\cite{Gilboa}) {\em For sufficiently large $n$,

(1) $ar(K_{n}, P_{3}\cup tP_{2})=(t-1)(n-\frac{t}{2})+1$ for $t\ge 2;$

(2) $ar(K_{n}, P_{4}\cup tP_{2})=t(n-\frac{t+1}{2})+1$ for $t\ge 1;$

(3) $ar(K_{n}, C_{3}\cup tP_{2})=t(n-\frac{t+1}{2})+1$ for $t\ge 1;$

(4) $ar(K_{n}, tP_{3})=(t-1)(n-\frac{t}{2})+1$ for $t\ge 1;$

(5) $ar(K_{n}, P_{k+1}\cup tP_{3})=(t+ \lfloor\frac{k}{2}\rfloor -1)(n-\frac{t+\lfloor k/ 2 \rfloor}{2})+1+ (k\mod 2)$ for $k\ge 3$ and $ t\ge 0;$

(6) $ar(K_{n}, P_{2}\cup tP_{3})=(t-1)(n-\frac{t}{2})+2$ for $t\ge 1;$

(7) $ar(K_{n}, kP_{2}\cup tP_{3})=(t+k-2)(n-\frac{t+k-1}{2})+1$ for $k\ge 2$ and $ t\ge 2.$
}

The Tur\'an number of star forests and linear forests are considered by Lidick\'y, Liu and Palmer \cite{Lidichy}.

\noindent{\bf Theorem 6.}(\cite{Lidichy}) {\em Let $F=\bigcup^{t}_{i=1}K_{1, p_{i}}$ be a star forest and $p_{1}\ge p_{2}\ge\cdots\ge p_{t}\ge 1$. For $n$ sufficiently large,
$$ex(n,F)=\max_{1\le i\le t}\bigg\{(i-1)n-\binom{i}{2}+\bigg\lfloor\frac{p_{i}-1}{2}(n-i+1)\bigg\rfloor\bigg\}.$$}

\noindent{\bf Theorem 7.}(\cite{Lidichy}) {\em Let $F=\bigcup^{k}_{i=1}P_{p_{i}}$ be a linear forest, where  $k\ge 2$ and $p_{i}\ge 2$ for $1\le i \le k$. If at least one $p_{i}$ is not $3$, then for $n$ sufficiently large, $$ex(n, F)=\left(\sum_{i=1}^{k}\left\lfloor\frac{p_{i}}{2}\right\rfloor-1\right)n-\binom{\sum_{i=1}^{k}\lfloor\frac{p_{i}}{2}\rfloor}{2}+c,$$ where $c=1$ if all $p_{i}$ are odd and $c=0$ otherwise.}

In Sections 2 and  3, we generalize Theorem 5 by considering the anti-Ramsey number of star forests and linear forests, respectively.

\noindent{\bf Theorem 8.} {\em Let $F=\bigcup^{t}_{i=1}K_{1, p_{i}}$ be a star forest, where $p_{1}\ge 3$, $p_{1}\ge p_{2}\ge\cdots\ge p_{t}\ge 1$. Let $s=\max\{i: p_{i}\ge 2, 1\le i\le t\}$. For  $n\ge 3t^{2}(p_{1}+1)^{2}$, we have
$$ar(K_n,F)=\max\bigg\{ \max_{1\le i\le s}\bigg\{(i-1)n-\binom{i}{2}+\bigg\lfloor\frac{p_{i}-2}{2}(n-i+1)\bigg\rfloor+1\bigg\},  (t-2)n-\binom{t-1}{2}+r\bigg\},$$
where $r = 1$ if $p_{t-1}=1$ and $r=2$ otherwise.}

\noindent{\bf Theorem 9.}  {\em Let $F=\bigcup^{k}_{i=1}P_{p_{i}}$ be a linear forest, where  $k\ge 2$ and $p_{i}\ge 2$ for $1\le i \le k$. We have
$$ar(K_{n}, F)= \left(\sum_{i=1}^{k}\left\lfloor\frac{p_{i}}{2}\right\rfloor-\epsilon\right)n+O(1),$$ where   $\epsilon=1$ if all $p_{i}$ are odd and $\epsilon=2$ otherwise.}

We get the approximate value of the anti-Ramsey number for linear forests by Theorem 9 and it would be interesting to determine the exact value. Bialostocki, Gilboa and Roditty \cite{Bialostocki} and independently Gorgol and  G\"orlich \cite{Gorgol} showed that $ar(K_{n}, 2P_{3})=\max\{n, 7\}$ for $n\ge 6$. Gorgol and  G\"orlich showed that $ar(K_{n}, 3P_{3})=2n-2$ for $n\ge 13.$ In Section 4, we will use Theorem 3 to compute the exact value of $ar(K_{n}, 2P_{4})$ for $n\ge 8$.

\noindent{\bf Theorem 10.} {\em For any $n\ge 8$, $ar(K_{n}, 2P_{4})=\max\{2n-2, 16\}.$}

Another motivation of this paper is the following conjecture of Gorgol and  G\"orlich \cite{Gorgol}:

{\em Let $G$ be a connected graph on $n_{0}\ge 3$ vertices and  $t\ge 1$, then for large $n$,  $$ar(K_{n}, tG)=(t-1)n-\binom{t}{2}+ar(K_{n-t+1},G)$$
if and only if $G$ is a tree.}

Statement (4) in Theorem 5(4) and Theorem 8 show that this conjecture is true for $P_{3}$ and $K_{1,p}$ ($p\ge 3$), respectively. However, from Theorem 9, some simple calculation shows that this conjecture fails for $P_{l}, l\ge 4$.

Actually, for an arbitrary tree $T_{k}$ with $k$ edges, it is difficult to determine the (approximate) value of $ar(K_{n}, T_{k})$. Jiang and West \cite{Jiang-1} showed that for $n\ge 2k$,
$$\bigg\lfloor\frac{k-2}{2}\bigg\rfloor\frac{n}{2}+O(1)\le ar(K_{n}, T_{k})\le (k-1)n.$$
The upper bound comes from the well-known bound of $ex(n, T_{k})\le (k-1)n.$
 Erd\H os and S\'os gave the following conjecture.

\noindent{\bf Conjecture 1.} $$ex(n, T_{k})\le \frac{k-1}{2}n.$$

If Conjecture 1 is true (Ajtai, Koml\'os, Simonovits, Szemer\'edi announced it for large $k$), then the upper bound of $ar(K_{n}, T_{k})$ can also be reduced to $\frac{k-1}{2}n.$  Also, Jiang and West \cite{Jiang-1} conjectured that:

\noindent{\bf Conjecture 2.} $$ar(n, T_{k})\le \frac{k-2}{2}n+O(1).$$

Notice that if $T_{k}$ is a star or a path of even length, then $ar(n, T_{k})= \frac{k-2}{2}n+O(1).$

The {\em double star} $S_{p,q}$, where $p\ge q\ge 1$, is the graph consisting of the union of two stars $K_{1,p}$ and $K_{1,q}$ together with an edge joining their centers.  In Section 5, we compute the anti-Ramsey number of double stars.

\noindent{\bf Theorem 11.} {\em For $p\ge 2, 1\le q\le p$ and $n\ge 6(p^{2}+2p),$ we have
\begin{align} \nonumber
ar(K_{n}, S_{p,q})=
 \begin{cases}
\lfloor\frac{(p-1)n}{2}\rfloor+1 ,\, 1\le q\le p-1;\\
\lfloor\frac{p(n-1)}{2}\rfloor+1,\, q=p.
 \end{cases}
 \end{align}}

Notice that if we take $T_{k}=S_{p,p-1}$, then we have $ar(K_{n}, T_{k})=\lfloor\frac{k-2}{2}\rfloor\frac{n}{2}+O(1).$

This paper is organized as follows. In Section 2, we give the proof of Theorem 8. The proof of Theorem 9 will be given in Section 3. The proof of Theorem 10 will be given in Section 4. The proof of Theorem 11 will be given in Section 5. Finally we will give a conjecture in Section 6.

{\bf Notation:}  Given an edge-coloring $c$ of $G$, we denote the color of an edge $uv$ by $c(uv)$. We denote the number of colors by $|c|$. For any $v\in V(G)$,  let $C(v):=\{c(vw)|\, w\in N_G(v)\}$ and $d_{c}(v):=|C(v)|.$  %For any color $a$, let $N_{c}(v; a):=\{w\in N_{G}(v)|\, c(vw)=a\}$.
Let $H$ be a subgraph of $G$. Denote $C(H)=\{c(uv)|\, uv\in E(H)\}$.
A color $a$ is  {\em stared} (at $x$) if all the edges with color $a$ induce a star $K_{1,r}$ (centered at the vertex $x$). We let $d^{c}(v)=|\{a\in C(v) |\, a $ is stared at $ v\}|$.
A {\em representing subgraph} in an edge-coloring of $K_n$ is a spanning subgraph
containing exactly one edge of each color. In the rest of this paper, we will use $V$ to denote the vertex set of $K_n$ for short.

%%%%%%%%%%%%%%%%%%%%%%%%%%%%%%%%%%%%%%%%%%%%%%
\vskip.2cm
\n{\large\bf 2. Star forests}
\vskip.2cm

In this Section, we use the idea of \cite{Gilboa} to prove Theorem 8.

\noindent{\bf Theorem 8.} {\em Let $F=\bigcup^{t}_{i=1}K_{1, p_{i}}$ be a star forest, where $p_{1}\ge 3, p_{1}\ge p_{2}\ge\cdots\ge p_{t}\ge 1$. Let $s=\max\{i: p_{i}\ge 2, 1\le i\le t\}$. For  $n\ge 3t^{2}(p_{1}+1)^{2}$,
$$ar(K_n,F)=\max\bigg\{ \max_{1\le i\le s}\bigg\{(i-1)n-\binom{i}{2}+\bigg\lfloor\frac{p_{i}-2}{2}(n-i+1)\bigg\rfloor+1\bigg\},  (t-2)n-\binom{t-1}{2}+r\bigg\},$$
where $r = 1$ if $p_{t-1}=1$ and $r=2$ otherwise.}

\begin{proof}
 For $1\le i \le s$, we color $K_{n}$ as follows. We color $K_{i-1}+\overline K_{n-i+1}$ rainbow and color $K_{n-i+1}$ with new $ar(K_{n-i+1}, K_{1,p_{i}})$ colors without producing a rainbow copy of $K_{1,p_{i}}$. In such way, we use exactly $(i-1)(n-i+1)+\binom{i-1}{2}+ar(K_{n-i+1}, K_{1,p_{i}})=(i-1)n-\binom{i}{2}+\big\lfloor\frac{p_{i}-2}{2}(n-i+1)\big\rfloor+1$  colors. Since any rainbow  $K_{1,p_{j}}$ ($ 1\le j\le i$) contains at least one vertex of $V(K_{i-1})$, we  do not obtain any rainbow  $\bigcup^{i}_{j=1}K_{1, p_{j}}$. Also, we do not obtain any rainbow  $F$.

For  another lower bound, it is enough to consider the case $t\ge 3$. If $p_{t-1}=1,$ then
$F\supset tK_{2}$. By Theorem 4, we have
$$ar(K_{n}, F)\ge ar(K_{n}, tK_{2})= ex(n, (t-1)K_{2})+1=(t-2)n-\binom{t-1}{2}+1.$$
If $p_{t-1}\ge 2$, we have $F\supset (t-1)K_{1,2}\cup P_{2}$. By  Theorem 5(6), we have
$$ar(K_{n}, F)\ge ar(K_{n}, (t-1)K_{1,2}\cup K_{2})=(t-2)n-\binom{t-1}{2}+2.$$

Now we consider the upper bound. If $t=1$, then the result holds obviously by Theorem 1. Suppose $t\ge 2$.
Let $$f(n, F)=\max\bigg\{ \max_{1\le i\le s}\bigg\{(i-1)n-\binom{i}{2}+\bigg\lfloor\frac{p_{i}-2}{2}(n-i+1)\bigg\rfloor+1\bigg\},  (t-2)n-\binom{t-1}{2}+r\bigg\}.$$
By Theorems 4 and 5 (statements (1), (4), (6), (7)), we have
 \begin{align} \tag{*}
 ar(K_{n}, kK_{1,2}\cup lK_{1,1})=
 \begin{cases}
 (k+l-2)n-\binom{k+l-1}{2}+1, \, k\ge 0, l\ge 2;\\
 (k-1)n-\binom{k}{2}+2, \, k\ge 1, l=1;\\
 (k-1)n-\binom{k}{2}+1,\,k\ge 1, l=0.
 \end{cases}
 \end{align}

 Let $c$ be any edge-coloring of $K_{n}$ using
$f(n, F)+1$ colors. We will find  a rainbow  $F$ by considering the following two cases.

{\bf Case 1.} There is some vertex $v_{0}\in V$ such that $d_{c}(v_{0})\ge \sum_{i=1}^{t}p_{i}+t$.

Choose some color $c_{0}$. Consider an edge-coloring $c'$ of $K_{n}-v_{0}$: $c'(e)=c(e)$ if $c(e) \notin C(v_{0})$; else $c'(e)=c_{0}$, where $e\in E(K_{n}-v_{0})$.

{\bf Claim 1.} There is a rainbow copy of $F-K_{1,p_{1}}$ in $K_{n}-v_{0}$ with respect to $c'$.

{\bf Proof of Claim 1.} We consider the following  two subcases.

%%%%%%%%%%%%%%%%%%%%%%%%%
{\bf Subcase 1.1} $p_{2}\le 2.$

In this subcase, we have $F-K_{1,p_{1}}=(s-1)K_{1,2}\cup (t-s)K_{1,1}$ and
$$f(n, F)=\max\bigg\{\bigg\lfloor\frac{p_{1}-2}{2}n\bigg\rfloor+1, (s-1)n-\binom{s}{2}+1, (t-2)n-\binom{t-1}{2}+r\bigg\}.$$
If $F-K_{1,p_{1}}=K_{1,1}$, then the result holds obviously. Suppose $F-K_{1,p_{1}}\ne K_{1,1}$. Then
\[\begin{split}
|c'| &\ge f(n, F)+1-(n-1)\\
&\ge \max\bigg\{(s-1)n-\binom{s}{2}+1, (t-2)n-\binom{t-1}{2}+r\bigg\}+1-(n-1)\\
&= \max\bigg\{(s-2)(n-1)-\binom{s-1}{2}+1, (t-3)(n-1)-\binom{t-2}{2}+r\bigg\}+1\\
&\ge ar(K_{n-1}, F-K_{1,p_{1}})+1,
\end{split}\]
by (*). Thus there is a rainbow  $F-K_{1,p_{1}}$ in $K_{n}-v_{0}$ with respect to $c'$.

{\bf Subcase 1.2} $p_{2}\ge 3.$

In this subcase, we have
$$|c'|\ge f(n, F)+1-(n-1)\ge f(n-1, F-K_{1,p_{1}})+1= ar(K_{n-1}, F-K_{1,p_{1}})+1.$$
By induction hypothesis,  there is a rainbow  $F-K_{1,p_{1}}$ in $K_{n}-v_{0}$ with respect to $c'$. \q

By Claim 1, there is a rainbow  $F-K_{1,p_{1}}$ in $K_{n}-v_{0}$ with respect to $c'$ (also respect to $c$).
Since $d_{c}(v_{0})\ge \sum_{i=1}^{t}p_{i}+t$, we are surely left with at least $p_{1}$ edges, say $v_{0}w_{1},v_{0}w_{2},\ldots,v_{0}w_{p_{1}}$, such that $w_{1},w_{2},\ldots,w_{p_{1}}\notin V(F-K_{1,p_{1}})$ and $c(v_{0}w_{1}),c(v_{0}w_{2}),\ldots,c(v_{0}w_{p_{1}})\notin \{c(e) | e\in E(F-K_{1,p_{1}})\} $. By adding such $p_{1}$ edges to $F-K_{1,p_{1}}$, we get a rainbow  $F$.

%%%%%%%%%%%%%%%%%%%%%%%%

%The vertex $v_{0}$ is the endpoint of at least $\sum_{i=1}^{t}p_{i}+t$ edges with distinct colors (with respect to $c$). There are at most $\sum_{i=2}^{t}p_{i}+t-1$ endpoints of those edges are in $V(H)$. Also, at most one of those edges has the same color (with respect to $c$) as an edge of $H$, since there is at most one edge of $H$ whose color is $c_{0}$  with respect to  $c'$. Therefore we are surely left with at least $p_{1}$ edges, say $v_{0}w_{1},v_{0}w_{2},\ldots,v_{0}w_{p_{1}}$, such that $w_{1},w_{2},\ldots,w_{p_{1}}\notin V(H)$ and $c(v_{0}w_{1}),c(v_{0}w_{2}),\ldots,c(v_{0}w_{p_{1}})\notin \{c(e) | e\in E(H)\} $. By adding such $p_{1}$ edges to $H$, we get a rainbow copy of $F$.

{\bf Case 2.} $d_{c}(v)\le  \sum_{i=1}^{t}p_{i}+t-1$ for all $v\in V$.

By induction hypothesis, $K_{n}$ clearly contains a rainbow  $F-K_{1,p_{t}}$.
Assume, by contradiction, that $K_{n} $ does not  contain a rainbow  $F$.
Let $G$ be a representing subgraph of $K_{n}$ such that $E(F-K_{1,p_{t}})\subseteq E(G)$. Then we have $|E(G)|=f(n,F)+1$ and $d_{G}(v)\le \sum_{i=1}^{t}p_{i}+t-1$  for all $ v \in V.$

Let $W=V-V(F-K_{1,p_{t}})$. Since $d_{G}(v)\le \sum_{i=1}^{t}p_{i}+t-1$ for all $v\in V(F-K_{1,p_{t}})$ and $G[W]$ does not contain a copy of $K_{1,p_{t}}$, we have
\[\begin{split}|E(G)| & \le |E_{G}(V(F-K_{1,p_{t}}), V)|+|E(G[W])|\\
& \le \big(\sum_{i=1}^{t-1}p_{i}+t-1\big)\big(\sum_{i=1}^{t}p_{i}+t-1\big)+\frac{(p_{t}-1)[n-(\sum_{i=1}^{t-1}p_{i}+t-1)]}{2} \\
& \le \big(\sum_{i=1}^{t}p_{i}+t-1\big)^{2} + \frac{p_{t}-1}{2}n.
\end{split}\]
We will finish the proof by considering the following two subcases.

{\bf Subcase 2.1} $p_{t}=1.$

Since $p_{1}\ge 3$ and $p_{t}=1$, we have
$$(\sum_{i=1}^{t}p_{i}+t-1)^{2} \ge |E(G)| =f(n, F)+1 \ge \bigg\lfloor\frac{p_{1}-2}{2}n\bigg\rfloor+2\ge \frac{n}{2}+1,$$
 a contradiction with $ n\ge 3t^{2}(p_{1}+1)^{2}.$

{\bf Subcase 2.2} $p_{t}\ge 2.$

In this subcase, we have  $s=t$ and
\[\begin{split}&\big(\sum_{i=1}^{t}p_{i}+t-1\big)^{2} + \frac{p_{t}-1}{2}n\ge |E(G)| =f(n, F)+1\\
   & \ge (t-1)n-\binom{t}{2}+\bigg\lfloor\frac{p_{t}-2}{2}(n-t+1)\bigg\rfloor+2 \\
   & \ge \frac{(2t+p_{t}-4)n}{2}-\frac{(t-1)(p_{t}+t-2)}{2}.
\end{split}\]
Thus we have $n\le \frac{2}{2t-3}[(\sum^{t}_{i=1}p_{i}+t-1)^{2}+(t-1)(p_{t}+t-2)],$
 a contradiction with $ n\ge 3t^{2}(p_{1}+1)^{2}.$
\end{proof}

%%%%%%%%%%%%%%%%%%%%%%%%%%%%%%%%%%%%%%%%%%%%%
\vskip.2cm
\n{\large\bf 3. Linear forests}
\vskip.2cm

First, we have the lower bound of the anti-Ramsey number for linear forests.

\noindent{\bf Proposition 1.}  {\em Let $F$ be a linear forest with components of order $p_{1}, p_{2}, \ldots, p_{k}$, where  $k\ge 2$ and $p_{i}\ge 2$ for $1\le i \le k$. Let $n\ge \sum_{i=1}^{k}p_{i}$. Then we have
$$ar(K_{n}, F)\ge\max\biggm\{ \binom{\sum_{i=1}^{k}p_{i}-2}{2}+1, sn-\binom{s+1}{2}+r\biggm\},$$
where  $s=\sum_{i=1}^{k}\lfloor\frac{p_{i}}{2}\rfloor-\epsilon$; $\epsilon=1$ if all $p_{i}$ are odd and $\epsilon=2$ otherwise; $r=2$ if exactly one $p_{i}$ is even and $r=1$ otherwise.}

\begin{proof}
For the first lower bound, we choose a subgraph $K_{\sum_{i=1}^{k}p_{i}-2}$ and color it rainbow. Then we use one extra color to color the remaining edges. In this way, we use exactly $\binom{\sum_{i=1}^{k}p_{i}-2}{2}+1$ colors and do not obtain a rainbow  $F$.

For the second lower bound, we  color $K_{s}+\overline{K}_{n-s}$ rainbow and color the edges of $K_{n-s}$ with $r$ new colors. Every copy of $F$ in $K_{n}$  have at least $(r+1)$ edges in $K_{n-s}$. In this way we do not obtain a rainbow  $F$ and use exactly $s(n-s)+\binom{s}{2}+r=sn-\binom{s+1}{2}+r$ colors.
\end{proof}

If all the components of the linear forest are even paths or odd paths, we can get the following corollary  from Theorems 2 and 7.

\noindent{\bf Corollary 1.} {\em Let $F$ be a linear forest with components of order $p_{1}, p_{2}, \ldots, p_{k}$, where  $k\ge 2$ and $p_{i}\ge 2$ for $1\le i \le k$ and $n$ sufficiently large. Let $s=\sum_{i=1}^{k}\lfloor\frac{p_{i}}{2}\rfloor-2$. If all $p_{i}$ are even, we have $$sn-\binom{s+1}{2}+1\le ar(K_n, F)\le sn-\binom{s+1}{2}+2.$$
If all $p_{i}$ are odd, we have $$ar(K_{n}, F)=(s+1)n-\binom{s+2}{2}+1.$$}

\begin{proof}The lower bound is due to Proposition 1.

For the upper bound, when all $p_{i}$ are even, by Theorem 2,  $$ar(K_{n}, F)\le ar(K_{n}, P_{\sum_{i=1}^{k}p_{i}})=sn-\binom{s+1}{2}+2.$$
When all $p_{i}$ are odd, if $p_{1}=p_{2}=\cdots=p_{k}=3$, by Theorem 5(4), $ar(K_{n}, F)=ar(K_{n}, kP_{3})=(k-1)n-\binom{k}{2}+1$. If at least one $p_{i}$ is not $3$, by Theorem 7, $$ar(K_{n}, F)\le ex(n, F)=(s+1)n-\binom{s+2}{2}+1.$$
\end{proof}

It is enough to consider the linear forests with at least one even path.

\noindent{\bf Theorem 12.} {\em Let $F$ be a linear forest with components of order $p_{1}, p_{2}, \ldots, p_{k}$, where  $k\ge 1$, $p_{i}\ge 2$ for $1\le i \le k$ and at least one $p_{i}$ is even. Then
$$ar(K_{n}, F)=\left(\sum_{i=1}^{k}\lfloor\frac{p_{i}}{2}\rfloor-2\right)n+O(1).$$}

\begin{proof}

By Proposition 1, we just need to show the upper bound.
We will use the idea of \cite{Simonovits} to prove it. The following results of Erd\H os and Gallai \cite{Erdos1} will be used in our proof.
$$\text{(a)} \,~~~~~~~~~~~~~~~~~~~~~~~~~ex(n, P_{r})\le \frac{r-2}{2}n;$$
$$\text{(b)}\,~~~~~~~ ex(n,\{C_{r+1}, C_{r+2},\cdots\} )\le \frac{r(n-1)}{2}.$$
Since we can regard the union of even paths as the subgraph of one long even path (see Corollary 1), we just need to prove the upper bound is correct for linear forest with exact one even path and some odd paths. So we assume that $F=P_{2s}\cup P_{2t_{1}+1}\cup P_{2t_{2}+1}\cup \cdots \cup P_{2t_{k}+1}$, where $s\ge 1, k\ge 1$ and $t_{k}\ge t_{k-1}\ge\cdots\ge t_{1}\ge1$. Let $s+t_{1}+t_{2}+\cdots+t_{k}=m$. In this proof, we just consider the case $s\ge 6$. The case $s<6$ can be proved by the similar arguments but need to distinguish more cases as in \cite{Simonovits}.

Consider an edge-coloring  of $K_{n}$ with $ar(K_{n}, F)$ colors such that there is no rainbow  $F$. First we take a rainbow path $P_{l}=u_{1}u_{2}\ldots u_{l}$ with maximum length. If $l\ge 2m+k$, we can get a rainbow $F$, a contradiction. By Proposition 1 and Theorem 2, $ar(K_{n}, F)\ge ar(K_{n}, P_{2m-1})+1$ for large $n$, which implies  there is a rainbow $P_{2m-1}$. Hence we  assume that $2m-1\le l\le2m+k-1.$ Take a representing subgraph $G$ of $K_{n}$ such that $P_{l}\subset G$. Then $|E(G)|=ar(K_n,F)$. We would partition $V\setminus V(P_{l})$ into three sets $U_{1}, U_{2}$ and $U_{3}$ as follows:

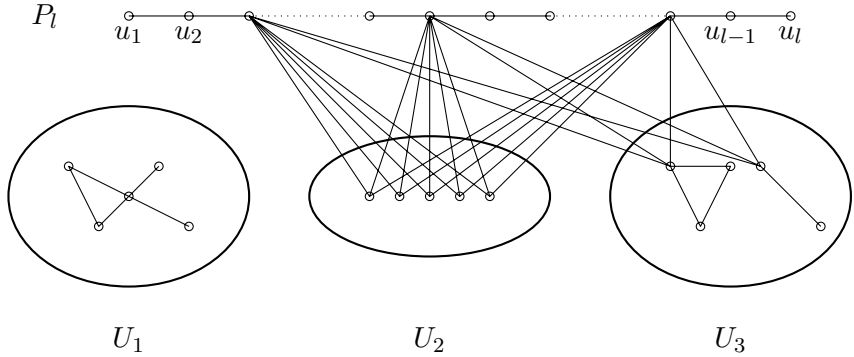
\begin{center}

\begin{tikzpicture}[scale=.8]

\draw[,thick] (0,0) arc(360:0:2cm and 1.5cm)  (5,0) arc(360:0:2cm and 1cm) (10,0) arc(360:0:2cm and 1.5cm)(-2,-2)node[below]{$U_{1}$}(3,-2)node[below]{$U_{2}$}(8,-2)node[below]{$U_{3}$};

\draw (-3,0.5)--(-2,0)--(-1.5,0.5) circle (2pt) (-3,0.5) circle (2pt)-- (-2.5,-0.5) circle (2pt)-- (-2,0) circle (2pt) --(-1,-0.5) circle (2pt)
(2,0) circle (2pt) (2.5,0) circle (2pt) (3,0) circle (2pt) (3.5,0) circle (2pt) (4,0) circle (2pt)
(7,0.5) circle (2pt)--(7.5,-0.5) circle (2pt)--(8,0.5) circle (2pt)--(7,0.5) (8.5,0.5) circle (2pt)--(9.5,-0.5) circle (2pt);

\draw (-3,3) node[left]{$P_{l}$} (-2,3)node[below]{$u_{1}$} circle (2pt)--(-1,3)node[below]{$u_{2}$}  circle (2pt)--(0,3) circle (2pt) (2,3) circle (2pt)--(3,3) circle (2pt)--(4,3) circle (2pt)--(5,3) circle (2pt) (7,3) circle (2pt)--(8,3)node[below]{$u_{l-1}$}  circle (2pt)--(9,3)node[below]{$u_{l}$}  circle (2pt);

\draw [rounded corners,dotted] (0,3) -- (2,3) (5,3) -- (7,3);

\draw (2,0)--(0,3) (2,0)--(3,3) (2,0)--(7,3)
(2.5,0)--(0,3) (2.5,0)--(3,3) (2.5,0)--(7,3)
(3,0)--(0,3) (3,0)--(3,3) (3,0)--(7,3)
(3.5,0)--(0,3) (3.5,0)--(3,3) (3.5,0)--(7,3)
(4,0)--(0,3) (4,0)--(3,3) (4,0)--(7,3)
(7,0.5)--(0,3) (7,0.5)--(3,3) (7,0.5)--(7,3)
(8.5,0.5)--(0,3) (8.5,0.5)--(3,3) (8.5,0.5)--(7,3);

\end{tikzpicture}

\end{center}

$U_{1}$ is the subset of vertices of $V\setminus V(P_{l})$ which are not jointed to $P_{l}$ at all: neither by edges nor by paths;

$U_{2}$ is the set of isolated vertices of  $V\setminus V(P_{l})$ which are jointed to $P_{l}$ by edges;

$U_{3}=V\setminus (V(P_{l})\cup U_{1} \cup U_{2}).$

{\bf Claim 1.}  $|E(G[U_{1}])|\le (m-2)|U_{1}|$.

{\bf Proof of Claim 1}
We  first prove that there is an  $P_{2s}\cup P_{2t_{1}+1}\cup \ldots \cup P_{2t_{k}+1}$ in $P_{2s+2t_{1}+\ldots+2t_{k}-1}\cup P_{2s+2t_{1}+\ldots+2t_{k}-2}$   by induction on $k$.  The base case $k=1$ is correct since $P_{2s}\cup P_{2t_{1}+1}\subset P_{2s+2t-1}\cup P_{2s+2t_{1}-2}$. Suppose the statement holds for $k-1$. We divide $P_{2s+2t_{1}+\ldots+2t_{k}-1}$ and $P_{2s+2t_{1}+\ldots+2t_{k}-2}$ respectively  into two parts $P_{2s+2t_{1}+\ldots+2t_{k-1}-1}, P_{2t_{k}}$ and $P_{2s+2t_{1}+\ldots+2t_{k-1}-2}, P_{2t_{k}}$.
We can find an  $P_{2s}\cup P_{2t_{1}+1}\cup \ldots \cup P_{2t_{k-1}+1}$ in $P_{2s+2t_{1}+\ldots+2t_{k-1}-1}\cup P_{2s+2t_{1}+\ldots+2t_{k-1}-2}$ by  induction hypothesis. Since $(2s+2t_{1}+\ldots+2t_{k-1}-1)+(2s+2t_{1}+\ldots+2t_{k-1}-2)> 2s +(2t_{1}+1)+\ldots+(2t_{k-1}+1)$,  there is  at least one vertex of either $P_{2s+2t_{1}+\ldots+2t_{k-1}-1}$ or  $P_{2s+2t_{1}+\ldots+2t_{k-1}-2}$ which is not used in  $P_{2s}\cup P_{2t_{1}+1}\cup \ldots \cup P_{2t_{k-1}+1}$. Hence, we can find an  $P_{2s}\cup P_{2t_{1}+1}\cup \ldots \cup P_{2t_{k}+1}$ in the original two long paths.

Since $G$ contains no $F$, $G[U_{1}]$ contains no $P_{2m-2}$  by the statement above.
 Thus  $|E(G[U_{1}])|\le (m-2)|U_{1}|$ by (a).\q

{\bf Claim 2.} $|\{v\in U_{2}: d_{G}(v)\ge m-1\}|\le (m+k)\binom{l-2}{m-1}$.

{\bf Proof of Claim 2}
It is obvious that  $N_{G}(v)\subset V(P_{l})\setminus\{u_{1}, u_{l}\}$ for all $v\in U_{2}$.
Suppose $|\{v\in U_{2}: d_{G}(v)\ge m-1\}|\ge (m+k)\binom{l-2}{m-1}+1$.  Note that there are $\binom{l-2}{m-1}$ subsets of order $m-1$ in $V(P_{l})\setminus\{u_{1}, u_{l}\}$. By  Pigeonhole Principle, there are $A\subseteq V(P_{l})\setminus\{u_{1}, u_{l}\}$ with $|A|=m-1$ and $B\subseteq U_{2}$ with $|B|=m+k+1$ such that $A\subseteq N_G(v)$ for any $v\in B$. So there is a rainbow $K_{m-1, m+k+1}$. By adding one edge whose endpoints are in the large part of $K_{m-1, m+k+1}$ to $K_{m-1, m+k+1}$, we can find a rainbow $F$, a contradiction.
\q

The following claim 3 is Lemma 1 in \cite{Simonovits}. We include the proof for the sake of completeness.

{\bf Claim 3.} $|E(G[U_{3}])|+|E_{G}(P_{l}, U_{3})|\le (m-2)|U_{3}|$.

{\bf Proof of Claim 3}
Take a component $H$ of $G[U_{3}]$ and  denote by $r$ the length of the longest cycle of $H$. If $H$ is a tree, let $r=2$. By (b), we have $|E(H)|\le \frac{r(|V(H)|-1)}{2}$. For any $v\in V(H)$, we can find an $P_{r}$ in $H$ starting from it. Hence, $u_{1}, \ldots, u_{r}$ and $u_{l-r+1}, \ldots, u_{l}$ cannot be joined to $v$. Otherwise, there is a rainbow $P_{l+1}$, a contradiction. For any three consecutive vertices $\{u_{i}, u_{i+1}, u_{i+2}\}$, there is no two independent edges in $E_{G}(\{u_{i}, u_{i+1}, u_{i+2}\}, V(H))$ by the maximality of $P_{l}$. Hence, we have
$$|E(H)|+|E_{G}(P_{l}, H)|\le \frac{r(|V(H)|-1)}{2}+\frac{l-2r+2}{3}|V(H)|\le\frac{l+1}{3}|V(H)|.$$
Adding all the components of $G[U_{3}]$ up, we get $|E(G[U_{3}])|+|E_{G}(P_{l}, U_{3})|\le \frac{l+1}{3}|U_{3}|$. Note that $l\le 2m+k-1$. Since $s\ge 6$, $m\ge s+k\ge 6+k$ which implies $|E(G[U_{3}])|+|E_{G}(P_{l}, U_{3})|\le (m-2)|U_{3}|$
\q

%%%%%%%%%%%%%%%

By Claims 1, 2 and 3, we have
\[\begin{split}ar(K_n,F) & =|E(G)|= |E(G[P_{l}])|+|E(G[U_{1}])|+|E_{G}(U_{2},P_{l})|+|E(G[U_{3}])|+|E_{G}(U_{3},P_{l})|\\
& \le \binom{l}{2}+ (m-2)|U_{1}|+ (l-2)(m+k)\binom{l-2}{m-1} +(m-2)|U_{2}|+(m-2)|U_{3}|\\
& \le  \binom{l}{2}+(l-2)(m+k)\binom{l-2}{m-1}+(m-2)(n-l)\\
&=(m-2)n+O(1).
\end{split}\]
\end{proof}

By  Corollary 1 and Theorem 12, we can get Theorem 9.

\vskip.2cm
\n{\large\bf 4. The exact value of $ar(K_{n}, 2P_{4})$}

\vskip.2cm
%%%%%%%%%%%%%%%%%%%%%%%%%%%%%%%%%%%%%%%
In this Section, we will prove Theorem 10. We  denote the complete graph on $n$ vertices minus one edge by $K_{n}^{-}$.
The following fact is trivial.

{\bf Fact 1.} Let $n\ge 8$. If there is  an edge coloring  of $K_{n}$ using $17$ colors such that there is a rainbow $K_{6}$ or $K_{6}^{-}$, then there is a rainbow $2P_{4}$.

 We first prove the following lemma.

\noindent{\bf Lemma 1.} {\em $ar(K_{8}, 2P_{4})=16.$}

\begin{proof} By Proposition 1 in Section 3, we just need to show the upper bound. Consider an $17$-edge-coloring $c$ of $K_{8}$. Suppose there is no rainbow  $2P_{4}$ in $K_{8}$. By Theorem 3, there must be a rainbow   $C_{k}\cup C_{l}$. Assume that $k\le l$. Then $k=3$. Let $T_{1}=C_3=x_1x_2x_3x_1$ and $T_{2}=C_l=y_1\ldots y_ly_1$. We choose a  representing subgraph $G$ such that $G\supset T_{1}\cup T_{2}$.  We just need to consider the following three cases.

{\bf Case 1.} $ l=5$.

 We claim that $c(xy)\in C(T_{1}\cup T_{2})$ for all $x\in V(T_{1})$ and $y\in V(T_{2})$. Otherwise,  say $c(x_{1}y_{1})\notin C(T_{1}\cup T_{2})$, then $x_{2}x_{3}x_{1}y_{1}\cup y_{2}y_{3}y_{4}y_{5}$ is a rainbow  $2P_{4}$, a contradiction. Then the total number of colors is at most $3+\binom{5}{2}=13<17$, a contradiction.

{\bf Case 2.} $l=4$.

Let  the remaining vertex be $z$. We have $c(x_{1}z)=c(x_{2}x_{3})$; otherwise there must be a rainbow  $2P_{4}$ in $T_{1}\cup T_{2}\cup\{x_{1}z\}.$ Similarly, we have $c(x_{2}z)=c(x_{1}x_{3})$ and $c(x_{3}z)=c(x_{1}x_{2})$. Now we have four rainbow $C_3$'s: $C_3^1=x_1x_2x_3x_1$, $C_3^2=zx_2x_3z$, $C_3^3=zx_1x_2z$ and $C_3^4=zx_1x_3x_1z$. If there are $u\in \{z,x_1,x_2,x_3\}$, say $u=z$,  and $1\le j\le 4$, say $j=1$, such that $c(zy_1)\not=c(zy_{2})$ and $c(zy_1),c(zy_{2})\notin C(C_3^1\cup T_2)$ then we have a rainbow $C_3^1$ and a rainbow $C_5=y_1zy_{2}y_3y_4 y_1$ and the situation is the same as Case 1. Hence we have $|E_{G}(u, T_{2})|\le2$ for any $u\in \{z,x_1,x_2,x_3\}$. Then  $17=|E(G)|\le3+4\times2+\binom{4}{2}=17$, which means $G[V(T_{2})]=K_4$ and $|E_{G}(u, T_{2})|=2$ for any $u\in \{z,x_1,x_2,x_3\}$. Thus we can find a rainbow  $2P_{4}$, a contradiction.

{\bf Case 3.} $l=3$.

Let  the remaining vertices be $z_{1}, z_{2}$. By Case 2 and $G$ containing no rainbow  $2P_4$, we have $|E_{G}(\{z_{1}, z_{2}\}, T_{1}\cup T_{2})|\le 2$ and then
 $ |E(G[V(T_{1}\cup T_{2})])|\ge 17-2-1=14.$ Thus
 $G[V(T_{1}\cup T_{2})]\cong K_{6}$ or $K_{6}^{-}$ and we can  get a rainbow  $2P_{4}$ by Fact 1, a contradiction.\end{proof}

Now we will complete the proof of Theorem 10.

\noindent{\bf Theorem 10.} {\em For any $n\ge 8$, $ar(K_{n}, 2P_{4})=\max\{2n-2, 16\}.$}

\begin{proof} By Proposition 1, we just need to show the upper bound. We  will prove it  by induction on $n$.

%%%%%%%%%%%%%%%   n=8
%\n{\bf 1. $ar(K_{8}, 2P_{4})=16.$}

 Consider an $(2n-1)$-edge-coloring of $K_{n}$ for $n\ge 9$. Suppose there is no rainbow  $2P_{4}$.   By Theorem 3, there is a rainbow  $C_{k}\cup C_{l}$. Assume that $k\le l$. Then $k=3$.  Let $T_{1}=C_3=x_1x_2x_3x_1$ and $T_{2}=C_l=y_1\ldots y_ly_1$.
 %For all $z\notin V(T_1)\cup V(T_2)$ when $4\le l\le n-4$, we have $c(x_{1}z)=c(x_{2}x_{3})$, $c(x_{2}z)=c(x_{1}x_{3})$ and $c(x_{3}z)=c(x_{1}x_{2})$; otherwise there is a rainbow copy of $2P_{4}$.
 We finish the proof by considering the following four cases.

{\bf Case 1.} $4\le l\le n-5.$

In this case, there are at least two vertices $z_{1}, z_{2} \notin V(T_{1}\cup T_{2})$.
Since there is no rainbow $2P_{4}$, we have $c(x_{1}z_{i})=c(x_{2}x_{3})$, $c(x_{2}z_{i})=c(x_{1}x_{3})$ and $c(x_{3}z_{i})=c(x_{1}x_{2})$ for $i=1, 2$. But we can have a rainbow  $2P_{4}$ whatever the color of $z_{1}z_{2}$ is, a contradiction.

{\bf Case 2.} $l= n-4$.

Let $V(K_{n})\setminus V(T_{1}\cup T_{2})=\{z\}.$
Since there is no rainbow  $2P_{4}$, we have that $c(x_{1}z)=c(x_{2}x_{3})$, $c(x_{2}z)=c(x_{1}x_{3})$, $c(x_{3}z)=c(x_{1}x_{2})$, $c(xy)\in C(T_{1}\cup T_{2})$ for all $x\in  V(T_{1})\cup \{z\}$ and $y\in V(T_{2})$. By Case 1, we can assume that  $c(y_{i}y_{j})\in C(T_{1}\cup T_{2})$ for all $y_{i}, y_{j}\in V(T_{2})$. Hence  the total number of colors is at most $3+n-4= n-1< 2n-1$, a contradiction.

{\bf Case 3.} $l= n-3$.

Since there is no rainbow  $2P_{4}$, we have that $c(xy)\in C(T_{1}\cup T_{2})$ for all $x\in  V(T_{1})$ and $y\in V(T_{2})$. By Case 1 and Case 2, we can assume that  $c(y_{i}y_{j})\in C(T_{1}\cup T_{2})$ for all $y_{i}, y_{j}\in V(T_{2})$. Hence  the total number of colors is at most $3+n-3= n< 2n-1$, a contradiction.

{\bf Case 4.} $l=3.$

In this case we will  consider the following two subcases. We choose a  representing subgraph $G$ such that $G\supset T_{1}\cup T_{2}$.

{\bf Subcase 4.1} $n=9.$

If there is a vertex $v_{0}$ such that $d^{c}(v_{0})=0$, then $K_{9}-v_{0}$ uses exactly $17$ colors  and there is a rainbow  $2P_{4}$ in $K_{9}-v_{0}$ by Lemma 1. Hence we have $d^{c}(v)\ge 1$, for all $v\in V(K_{9})$. Thus $\delta(G)\ge 1$.

Let   $Z=V(K_{9})\setminus V(T_{1}\cup T_{2})=\{z_{1}, z_{2}, z_{3}\}.$ We claim that $|E_G(z_j,T_s)|\le 1$ for any $1\le j\le 3$ and $1\le s\le2$. Otherwise $G$ contains an $C_3\cup C_4$ and the situation is the same as Case 1. If there are $z_i\in Z$ and $s\in \{1, 2\}$ such that $|E_G(z_i,T_s)|= 1$, then $|E_G(z,T_t)|= 0$ for any $z\in Z\setminus\{z_i\}$ and $t\in \{1,2\}\setminus\{s\}$ by $G$ having no  $2P_4$. If $E_{G}(Z, T_{1}\cup T_{2})=\0$, then $|E(G[V(T_{1}\cup T_{2})])|= 17-|E(G[Z])|\ge 17-3=14$ which implies $G[V(T_{1}\cup T_{2})]\cong K_{6}$ or $K_{6}^{-}$ and  there is a rainbow $2P_{4}$ by Fact 1, a contradiction. Hence we have $E_{G}(Z, T_{1}\cup T_{2})\ne\0$. Assume $z_1x_1\in E(G)$ and we will consider the following two subcases.

{\bf Subcase 4.1.1} $E_{G}(Z,  T_{2})=\0$.

Recall $z_1x_1\in E(G)$ and $|E_G(z_j,T_1)|\le 1$ for any $1\le j\le 3$. If there is $z\in \{z_2,z_3\}$ such that $z_1z\in E(G)$, then $E_{G}(\{x_{2},x_{3}\}, T_{2})$ $=\0$ which implies $|E(G)|\le 15$, a contradiction. Hence we have $|E(G[Z])|\le 1$ and then $|E(G[V(T_{1}\cup T_{2})])|\ge 17-3=14.$ Thus $G[V(T_{1}\cup T_{2})]\cong K_{6}$ or $K_{6}^{-}$ and  there is a rainbow  $2P_{4}$ by Fact 1, a contradiction.

{\bf Subcase 4.1.2}  $E_{G}(Z,  T_{2})\ne\0$.

In this case, we have $|E_{G}(Z,  T_{1}\cup T_{2})|=|E_{G}(z_{1},  T_1\cup T_{2})|=2$ and assume $z_1y_1\in E(G)$. If $E_{G}(z_{1}, \{z_{2}, z_{3}\})\ne \0, $  say $z_{1}z_{2} \in E(G)$, then $E_{G}(T_{1}, T_{2})\subseteq \{x_{1}y_{1}\}$. Hence $|E(G)|\le 12$, a contradiction. So we have $E_{G}(z_{1}, \{z_{2}, z_{3}\})= \0,$ which implies $z_{2}z_{3} \in E(G)$ by $\delta(G)\ge 1.$ Thus
 $|E(G[V(T_{1}\cup T_{2})])|= 17-3=14$
 and $G[V(T_{1}\cup T_{2})]\cong K_{6}^{-}$. We can get a rainbow  $2P_{4}$ by Fact 1, a contradiction.

{\bf Subcase 4.2} $n\ge 10.$

If there is a vertex $v_{0}$ such that $d^{c}(v_{0})\le 2$, then $K_{n}-v_{0}$ uses  $2n-1-d^{c}(v_{0})\ge 2n-1-2=2(n-1)-1$ colors and we get that $K_{n}-v_{0}$ contains a rainbow $2P_{4}$ by  induction hypothesis, a contradiction. So we assume that $d^{c}(v)\ge 3$ for all $v\in V(K_{n})$. Thus $\delta(G)\ge 3$.

Let  $Z=V(K_{n})\setminus V(T_{1}\cup T_{2})=\{z_{1}, z_{2},\ldots,  z_{n-6}\}.$ Since $G$ contains no  rainbow  $2P_{4}$, we have that if there are $z\in Z$ and $s\in \{1,2\}$ such that $E_{G}(z, T_{s})\ne \0$, then $E_{G}(z', T_{t})=\0$ for any $z'\in Z\setminus\{z\}$ and $t\in \{1,2\}\setminus\{s\}$.

{\bf Subcase 4.2.1} $E_{G}(Z, T_{1})\ne \0$ and $E_{G}(Z,T_{2})\ne \0.$

In this case, there is exactly one vertex in $Z$, say $z_{1}$, such that $E_{G}(z_{1}, T_{1}), E_{G}(z_{1},  T_{2})\ne \0$ and $E_{G}(Z\setminus\{z_{1}\}, T_{1}\cup T_{2})=\0$. Since $d_{G}(v)\ge 3$, $G[T_{1}\cup T_{2}\cup \{z_{1}\}]$ contains a cycle with length at least $4$ and  $G[Z\setminus\{z_{1}\}]$  contains a cycle of length at least 3. Then we can have a contradiction by Cases 1 to 3.

{\bf Subcase 4.2.2} $E_{G}(Z, T_{1})=\0$ or $E_{G}(Z, T_{2})=\0.$

Assume that $E_{G}(Z, T_{1})=\0.$ By Case 1, we can assume that $|E_{G}(z, T_{2})|\le 1$ for any $z\in Z$. Then $\delta(G[Z])\ge 2$ and there is a cycle in $G[Z]$. Since $d_{G}(x_{i})\ge 3$, we have $|E_{G}(x_{i}, T_{2})|\ge 1$ for any $1\le i\le 3$ and there is a cycle with length at least $4$ in $G[T_{1}\cup T_{2}]$. Then we can have a contradiction by Cases 1 to 3.
\end{proof}
%%%%%%%%%%%%%%%%%%%%%%%%%%%%%%%%%%%%%%%%%%%%%%%
\vskip.2cm

\n{\large\bf 5. Double stars}

\vskip.2cm

The following Lemma 2 is an extension of Theorem 1. The idea of the proof is the same as the idea used in \cite{Jiang}. We include the proof for the sake of completeness.

\noindent{\bf Lemma 2.} {\em Let $V(K_{n})=\{v_{1}, v_{2}, \ldots, v_{n}\}$. For $1\le p_{1}\le p_{2}\le\cdots\le p_{n}\le \frac{n}{3}$,  the maximum number of colors of an edge-coloring of $K_{n}$ such that $d_{c}(v_{i})\le p_{i}$ for any $ 1\le i\le n $ is at most $\frac{\sum^{n}_{i=1}p_{i} -n}{2}+1$.}

\begin{proof}

Let $c$ be an edge-coloring of $K_{n}$ such that $d_{c}(v_{i})\le p_{i}$ for any $ 1\le i\le n$. Let $G$ be a representing subgraph of $K_{n}$. Then $d_{G}(v_{i})\le d_{c}(v_{i})\le p_{i}$ for $ 1\le i\le n$ and $|E(G)|=|c|$.
Write $V=V(G)=V(K_{n})$. Let $S=\{v_{i}\in V: d_{G}(v_{i})=p_{i}\}$. If $S=\emptyset $, then $d_{G}(v_{i})\le p_{i}-1$ for $ 1\le i\le n$ and hence $$|c|=\frac{\sum^{n}_{i=1}d_{G}(v_{i})}{2}\le\frac{\sum^{n}_{i=1}(p_{i}-1)}{2}=\frac{\sum^{n}_{i=1}p_{i} -n}{2}$$and we are done.
Hence we may assume that $S\not=\emptyset$. For $v\in V$, let $C_{G}(v)$ be the set of colors used on the edges incident to $v$ in $G.$ Clearly, we have $C_{G}(v)\subset C(v) $ for all $v\in V$ and $C_{G}(v)=C(v)$ for $v\in S$. Particularly, $\{c(uv)\}= C(u)\cap C(v)= C_{G}(u)\cap C_{G}(v)$ for $u, v\in S$.

{\bf Claim 1.} $G[S]$ is a clique.

{\bf Proof of Claim 1.} Let $u, v\in S$. Since $c(uv)\in C(u)\cap C(v)= C_{G}(u)\cap C_{G}(v)$, $uv\in E(G)$. Hence $G[S]$ is a clique.\q

{\bf Claim 2.} Let $u, v\in S$ and $w\in V-S$. If $c(uw)=c(vw)$, then $c(uw)=c(vw)=c(uv)$.

{\bf Proof of Claim 2.} Since $c(uw)=c(vw)\in C(u)\cap C(v)= C_{G}(u)\cap C_{G}(v)=\{c(uv)\}$, we have $c(uw)=c(vw)=c(uv)$.\q

{\bf Claim 3.} Let $u\in S, v\notin S$ and $uv\notin E(G)$. Then $c(uv)\notin C_{G}(v)$.

{\bf Proof of Claim 3.} Suppose $c(uv)\in C_{G}(v)$. Since $u\in S$, we have $c(uv)\in C_{G}(u)$. Hence $c(uv)\in C_{G}(u)\cap C_{G}(v)$ which implies $uv\in E(G)$, a contradiction.\q

{\bf Claim 4.} For all $v_{i}\notin S$, $d_{G}(v_{i})\le p_{i}-\frac{|S\setminus N_{G}(v_{i})|}{2}$.

{\bf Proof of Claim 4.} Let  $S\setminus N_{G}(v_{i})=\{u_{1}, u_{2},\ldots, u_{k}\}.$ By Claim 3, $c(u_{j}v_{i})\notin C_{G}(v_{i})$ for all $1\le j\le k$. Furthermore, by Claim 2, if $c(u_{j}v_{i})=c(u_{l}v_{i})$, then $c(u_{j}v_{i})=c(u_{l}v_{i})=c(u_{j}u_{l})$. This implies that no three  edges in the set $\{u_{1}v_{i}, u_{2}v_{i}, \ldots, u_{k}v_{i}\}$  have the same color; otherwise, by Claim 1, $G[S]$ would contain a monochromatic triangle, a contradiction. Hence at least $\frac{k}{2}$ distinct colors are used on the edges $v_{i}u_{1}, v_{i}u_{2}, \cdots, v_{i}u_{k}$, and those colors are not  in $C_{G}(v_{i})$. So $|C(v_{i})|\ge |C_{G}(v_{i})|+\frac{k}{2}= d_{G}(v_{i})+\frac{k}{2}$.   Since $|C(v_{i})|=d_{c}(v_{i})\le p_{i}$, we have $d_{G}(v_{i})+\frac{k}{2}\le p_{i}$, which yields $d_{G}(v_{i})\le p_{i}-\frac{k}{2}=  p_{i}-\frac{|S\setminus N_{G}(v_{i})|}{2}.$\q

By Claim 4, we have
\[\begin{split}\sum_{i=1}^{n}d_{G}(v_{i})&=\sum_{v_{i}\in S}d_{G}(v_{i})+\sum_{v_{i}\notin S}d_{G}(v_{i})\\
&\le \sum_{v_{i}\in S}p_{i}+\sum_{v_{i}\notin S}\left(p_{i}-\frac{|S\setminus N_{G}(v_{i})|}{2}\right)\\
&=\sum_{i=1}^{n}p_{i}-\sum_{v_{i}\notin S}\frac{|S\setminus N_{G}(v_{i})|}{2}.\end{split}\]
Notice that $\sum_{v_{i}\notin S}|S\setminus N_{G}(v_{i})|$ counts exactly the number of non-edges in $G$ between $S$ and $V\setminus S$. We have $\sum_{v_{i}\notin S}|S\setminus N_{G}(v_{i})|=\sum_{v_{i}\in S}(n-1-p_{i}).$ Hence,
$$\sum_{i=1}^{n}d_{G}(v_{i})\le \sum_{i=1}^{n}p_{i}-n+n-\frac{\sum_{v_{i}\in S}(n-1-p_{i})}{2}.$$

On the other hand, for $v_{i}\notin S, d_{G}(v_{i})\le p_{i}-1$. Hence, we have
$$\sum_{i=1}^{n}d_{G}(v_{i})\le \sum_{v_{i}\in S}p_{i}+\sum_{v_{i}\notin S}(p_{i}-1)\le \sum_{i=1}^{n}p_{i}-n+|S|.$$

We have $$\sum_{i=1}^{n}d_{G}(v_{i})\le \sum_{i=1}^{n}p_{i}-n+\min\left\{|S|, n-\frac{\sum_{v_{i}\in S}(n-1-p_{i})}{2}\right\}.$$
Since $n\ge 3p_{n}$, we have $\min\{|S|, n-\frac{\sum_{v_{i}\in S}(n-1-p_{i})}{2}\}\le 2$. Therefore, $$|c|=\frac{\sum_{i=1}^{n}d_{G}(v_{i})}{2}\le \frac{\sum^{n}_{i=1}p_{i} -n}{2}+1.$$
\end{proof}

%%%%%%%%%%%%%%%%%%%%%%%%%%%%%%%%%%%%%%%%%

Now  we will  prove Theorem 11.

\noindent{\bf Theorem 11.} For $p\ge 2, 1\le q\le p$ and $n\ge 6(p^{2}+2p),$
\begin{align} \nonumber
ar(K_{n}, S_{p,q})=
 \begin{cases}
\lfloor\frac{(p-1)n}{2}\rfloor+1 ,\, 1\le q\le p-1;\\
\lfloor\frac{p(n-1)}{2}\rfloor+1,\, q=p.
 \end{cases}
 \end{align}

 \begin{proof}
 If $1\le q\le p-1$, we have $$ar(K_{n}, S_{p,q})\ge ar(K_{n}, K_{1,p+1})=\left\lfloor\frac{(p-1)n}{2}\right\rfloor+1 .$$ When $q=p,$ we color the edges of $K_{n}$ as follows. We color  $K_{1,n-1}$ rainbow and color $K_{n-1}$ with new $ar(K_{n-1}, K_{1,p})$ colors without producing a rainbow  $K_{1,p}$. In such way, we use exactly $$n-1+ar(K_{n-1}, K_{1,p})=n-1+\left\lfloor\frac{(p-2)(n-1)}{2}\right\rfloor+1=\left\lfloor\frac{p(n-1)}{2}\right\rfloor+1$$ colors and do not obtain a rainbow  $S_{p,p}$.

Now we consider the upper bound.  Let
 \begin{align} g(n, p, q)=
 \begin{cases}
\lfloor\frac{(p-1)n}{2}\rfloor+1 ,\, 1\le q\le p-1,\\
\lfloor\frac{p(n-1)}{2}\rfloor+1,\, q=p.
 \end{cases}
 \end{align}
Let $c$ be any edge-coloring of $K_{n}$ using $g(n,p,q)+1$ colors. We will find a rainbow $S_{p,q}$ by considering the following two cases.

{\bf Case 1.}  There is some vertex $v_{0}$ such that $d_{c}(v_{0})\ge p+2q+1.$

Let $U=\{ v\in V\setminus \{v_{0}\}: $ the color $c(v_{0}v)$ is stared at $ v_0 \}.$
Consider the induced edge-coloring $c'$ of $c$ on $K_{n}-v_{0}$. We have $$|c'|=|C(K_{n}-v_{0})|=g(n,p,q)+1-d^{c}(v_{0})\ge g(n,p,q)+1-|U|.$$

{\bf Claim 1.} $K_{n}-v_{0}$ contains a rainbow  $K_{1,q}$ with center in $U$ or  a rainbow  $K_{1,q+1}$ with center in $V(K_{n}-v_{0})\setminus U$ with respect to $c'$ (also $c$).

{\bf Proof of Claim 1.} Suppose $d_{c'}(u)\le q-1$ for all $u\in U$ and $d_{c'}(v)\le q$ for all $v\in V(K_{n}-v_{0})\setminus U$, by Lemma 2, we have
$$|c'|\le\frac{|U|(q-1)+(n-1-|U|)q-(n-1)}{2}+1.$$
While $|c'|\ge g(n, p, q)+1-|U|,$
we can get $|U|\ge n$, a contradiction.\q

By  Claim 1, we get a rainbow  $K_{1, q}$ in $K_{n}-v_{0}$ with respect to $c$ whose center is $u_{0}$  such that the color $c(v_{0}u_{0})$ does not  present in  $c(K_{1, q})$.  The vertex $v_{0}$ is the endpoint of at least $p+2q+1$ edges with distinct colors (with respect to $c$). Since $d_{c}(v_{0})\ge p+2q+1,$  there are at least $p$ edges, say $v_{0}v_{1}, v_{0}v_{2}, \ldots, v_{0}v_{p}$, such that $v_{i}\notin V(K_{1, q})$ and $c(v_{0}v_{i})\notin C(K_{1, q})\cup \{c(v_{0}u_{0})\}$ for all $i= 1,\ldots, p$. So we can get a rainbow  $S_{p,q}$.

{\bf Case 2.}  $d_{c}(v)\le p+2q$ for all $v\in V$.

Since $g(n, p, q)+1\ge ar(K_{n}, K_{1,p+1})+1$, we can find a rainbow  $K_{1,p+1}$. Choose some color $c_{0}$ which is not in $C(K_{1,p+1})$. Consider an edge-coloring $c''$ of $K_n-V(K_{1,p+1})$: for $e\in E(K_n-V(K_{1,p+1}))$, if $c(e)\notin C(K_{1,p+1})$, then $c''(e)=c(e)$; else $c''(e)=c_{0}$. Then
$$|c'|\ge g(n, p, q)+1-(p+2q)(p+2)\ge ar(K_{n-p-2}, K_{1,q+1})+1.$$
So there is a rainbow $K_{1,q+1}$ in $K_{n}-K_{1,p+1}$ with respect to $c''$ (also $c$).  $K_{1,q+1} $ has at most one edge $e$ with color $c_{0}$. By joining the centers $K_{1,p+1}$ and $K_{1,q+1}$ and deleting one edge of $K_{1,p+1}$ and $K_{1,q+1}$ respectively,  we can find a rainbow  $S_{p,q}$.
\end{proof}

\vskip.2cm

\n{\large\bf 6. Open problems}

\vskip.2cm

A {\em spider} is a tree with at most one vertex of degree more than $2$, called the center of the spider (if no vertex of degree more than two, then any vertex can be the center). A {\em leg} of a spider is a path from the center to  a vertex of degree 1. Thus, a star with $p$ edges is a spider of $p$ legs, each of length $1$, and a path is a spider of $1$ or $2$ legs.

%A {\em subdivision} of an edge $uv$ is achieved by removing the edge  and replacing it with a new vertex $w$ and the edges $uw$ and $wv.$ %If each pendant edge of the star $K_{1,p}$ is subdivided one time, the resulting graph will be defined as a {\em $1$-extended star} and denoted by $S^{1}_{p}$.
%If each pendant edge of the star $K_{1,p}$ is subdivided $t$ times, the resulting graph will be defined as an {\em $t$-extended star} and denoted by $S^{t}_{p}$. Notice that $S_{p}^{t}$ is a spider of $p$ legs, each of length $t+1$.

The number of edges in a maximum matching of a graph $G$ is called the {\em matching number} of $G$ and denoted by $\nu(G)$.

\noindent{\bf Observation 1.}  %Let $p\ge 3$ and $T$ be a spider of $p$ legs, each of length at least $2$. If $T$ has exactly one leg with  length even, then $$\min\{\nu(T-e_{1}-e_{2}): e_{1}, e_{2}\in E(T)\}=\min\{\nu(T-e): e\in E(T)\}.$$
{\em Let $p\ge 2$ and $T$ be a spider of $p$ legs, each of length at least $2$. We have  $$\min\{\nu(T-e_{1}-e_{2}): e_{1}, e_{2}\in E(T)\}\le \min\{\nu(T-e): e\in E(T)\},$$
the equality holds if and only if $T$ has exactly one leg with length even.}

Let $p\ge 2$ and $T$ be a spider of $p$ legs, each of length at least $2$. Let $\beta(T)=\min\{\nu(T-e): e\in E(T)\}$. %Note that if all the lengths of  the legs of $T$ are even (resp. odd), then $\beta(T)=\frac{e(T)}{2}$ (resp. $\beta(T)=\frac{e(T)-p+2}{2}$), where $e(T)=|E(T)|$.

\noindent{\bf Proposition 2.} {\em Let $p\ge 2$ and $T$ be a spider of $p$ legs, each of length at least $2$.  We have
$$ar(K_{n}, T)\ge (\beta(T)-1)n-\binom{\beta(T)}{2}+r,$$
where $r=2$ if there is exactly one leg of $T$ with length even and $r=1$ otherwise.}

\begin{proof}

Let $\beta=\beta(T)$. We take an $K_{\beta-1}+\overline{K}_{n-\beta+1}$ and color it rainbow, and use $r$ extra colors for all the remaining edges. Suppose there is a rainbow $T$ in this coloring. Then
$T-e$ contains a matching of size $\beta$ for any $e\in E(T)$ (or $T-e_{1}-e_{2}$ contains a matching of size $\beta$ for any $e_{1}, e_{2}\in E(T)$ if $T$ has exactly one leg with  length even). But $K_{\beta-1}+\overline{K}_{n-\beta+1}$ does not contain a matching of size $\beta$, a contradiction.
\end{proof}

 If we regard $P_{k+1}$ as a spider of $2$ legs, the lower bound of Proposition 2 is sharp for  $p=2$ and large $n$ by Theorem 2.
We  conjecture that the lower bound is sharp for $p\ge 3$ and large $n$.

\noindent{\bf Conjecture 3.}
 {\em Let $p\ge 3$ and $T$ be a spider of $p$ legs, each of length at least $2$.  For large $n$,
$$ar(K_{n}, T)= (\beta(T)-1)n-\binom{\beta(T)}{2}+r,$$
where  $r=2$ if there is exactly one leg of $T$ with length even and $r=1$ otherwise.}

%%%%%%%%%%%%%%%%%%%%%%%%%%%%%%%%%%%%

\end{document}